\input graphicx
\input amssym.def
\input amssym
\magnification=1000
\baselineskip = 0.18truein
\lineskiplimit = 0.01truein
\lineskip = 0.01truein
\vsize = 8.7truein
\voffset = 0.1truein
\parskip = 0.10truein
\parindent = 0.3truein
\settabs 12 \columns
\hsize = 6.0truein
\hoffset = 0.1truein

\setbox\strutbox=\hbox{%
\vrule height .708\baselineskip
depth .292\baselineskip
width 0pt}
\font\caps=cmcsc10

\font\bigtenrm=cmr10 at 14pt

\def\sqr#1#2{{\vcenter{\vbox{\hrule height.#2pt
\hbox{\vrule width.#2pt height#1pt \kern#1pt
\vrule width.#2pt}
\hrule height.#2pt}}}}
\def\square{\mathchoice\sqr46\sqr46\sqr{3.1}6\sqr{2.3}4}

\centerline{\bigtenrm EVERY KNOT HAS CHARACTERISING SLOPES}
\tenrm
\vskip 14pt
\centerline{MARC LACKENBY\footnote{\dag}{Partially supported by EPSRC grant EP/R005125/1}}
\vskip 18pt
\centerline{\caps 1. Introduction}
\vskip 6pt

Any closed orientable 3-manifold is obtained by Dehn surgery on a link in the 3-sphere [19, 32]. This
surgery description of the manifold is known to be highly non-unique [15]. Nevertheless, if one
restricts to knots, rather than links, and one considers only certain surgery slopes, then
uniqueness results are a realistic goal. One way of formalising this is via the notion of
a `characterising slope' which is defined as follows.

For a knot $K$ in $S^3$, we let $S^3_K(p/q)$ denote the 3-manifold obtained by $p/q$ Dehn surgery on $K$.
A slope $p/q$ is a {\sl characterising slope} for $K$ if whenever
there is an orientation-preserving homeomorphism between $S^3_K(p/q)$ and $S^3_{K'}(p/q)$
for some knot $K'$ in $S^3$, then there is an orientation-preserving homeomorphism of
$S^3$ taking $K$ to $K'$. (Recall that slopes on a boundary torus of the exterior of
a link in the 3-sphere are parametrised by fractions $p/q \in {\Bbb Q} \cup \{ \infty \}$, where 
$1/0 = \infty$ represents the meridian and $0/1 = 0$ represents the longitude of the relevant link component.
Throughout this paper the integers $p$ and $q$ will be assumed to be coprime.)

There has been considerable interest in characterising slopes over many years. A conjecture of Gordon [10],
proved by Kronheimer, Mrowka, Ozsv\'ath and Szab\'o [16], was that every slope is a characterising slope 
for the unknot. Using Heegaard Floer
homology, this was extended to the figure-eight knot and the trefoil by Ozsv\'ath and Szab\'o [24]. 
Then Ni and Zhang [23] showed
that certain slopes are characterising slopes for the torus knot $T_{5,2}$. This was developed
by McCoy [20], who showed that every torus knot has characterising slopes; in fact, all but finitely many
non-integral slopes are characterising for a torus knot. Ni and Zhang conjectured
that, for any hyperbolic knot $K$, all but finitely many slopes are characterising for $K$.
However, this was disproved by Baker and Motegi [2], who gave examples of hyperbolic knots ($8_6$ being the simplest)
for which no integer slope is characterising. Baker and Motegi conjectured that, nevertheless,
every knot has infinitely many characterising slopes. The main result of this paper is a
proof of their conjecture.

\noindent {\bf Theorem 1.1.} {\sl Let $K$ be a knot in $S^3$. Then any slope $p/q$ is characterising
for $K$, provided $|p| \leq |q|$ and $|q|$ is sufficiently large.}

Baker and Motegi also conjectured that for any hyperbolic knot $K$ in $S^3$, all but finitely many of its
non-integral slopes are characterising. We cannot prove this, but we do achieve the following
result.

\noindent {\bf Theorem 1.2.} {\sl Let $K$ be a hyperbolic knot in $S^3$. Then any slope $p/q$ is characterising
for $K$ provided that $|q|$ is sufficiently large.}

In the proof of both Theorem 1.1 and 1.2, one must consider not just the knot $K$ but also another knot $K'$.
The argument works most efficiently when both $K$ and $K'$ are hyperbolic, in which case we have
the following stronger version of Theorem 1.2.

\noindent {\bf Theorem 1.3.} {\sl Let $K$ be a hyperbolic knot in $S^3$. Then there is a constant $C(K)$
with the following property. If there is an orientation-preserving homeomorphism between
$S^3_K(p/q)$ and $S^3_{K'}(p'/q')$ for some hyperbolic knot $K'$ in $S^3$ and $|q'| \geq C(K)$, 
then there is an orientation-preserving homeomorphism of $S^3$ taking $K$ to $K'$
and $p/q = p'/q'$.}

A possible approach to Theorem 1.3 is to show that $S^3_K(p/q)$ `remembers' the Dehn surgery in the following sense.
By Thurston's hyperbolic Dehn surgery theorem [30], $S^3_K(p/q)$ is hyperbolic provided $|p|+|q|$ is sufficiently
large and moreover, the core of the surgery solid torus is the shortest geodesic in $S^3_K(p/q)$.
Suppose the same is true of $S^3_{K'}(p'/q')$ for some hyperbolic knot $K'$ in $S^3$. By Mostow rigidity [21], the homeomorphism
from $S^3_K(p/q)$ to $S^3_{K'}(p'/q')$ is homotopic to an isometry, which must therefore take
the shortest geodesic to the shortest geodesic. This therefore restricts to an orientation-preserving
homeomorphism from $S^3 - {\rm int}(N(K))$ to $S^3 - {\rm int}(N(K'))$ taking the slope
$p/q$ to the slope $p'/q'$. Using a theorem of Gordon and Luecke [11] (which is essentially that
knots in the 3-sphere are determined by their complements), this homeomorphism takes
the meridional slope on $\partial N(K)$ to the meridional slope on $\partial N(K')$, and hence extends to an orientation-preserving
homeomorphism from $S^3$ to $S^3$ taking $K$ to $K'$. This sort of argument has been used before, for example
by Bleiler, Hodgson and Weeks [5] and by Rieck and Yamashita [28]. They analysed the related concept of cosmetic surgeries,
where there is an orientation-preserving homeomorphism between $S^3_K(p/q)$ and $S^3_K(p'/q')$ 
but $p/q \not= p'/q'$ (and so the knot type is fixed but the slope can vary).

There is a significant difficulty with the line of argument presented in the previous paragraph. The core of the surgery solid torus is
indeed the shortest geodesic in $S^3_K(p/q)$ when $|p| + |q|$ is greater than some constant, but this constant 
depends on $K$. We are also considering surgery along $K'$ and here the relevant constant
depends on $K'$. Thus, we cannot immediately deduce that any slope $p/q$ is a characterising
slope for $K$. Indeed, the examples due to Baker and Motegi [2] of hyperbolic knots with infinitely
many non-characterising slopes highlight that this argument cannot be made to work in this form.
Nevertheless, by delving deeper into the theory of hyperbolic 3-manifolds and Dehn surgery, it is possible
to construct a proof of Theorem 1.3. The proof is given in Section 3, after the relevant background material on Dehn surgery has been collated
in Section 2.

To prove Theorem 1.2, we must also consider the possibility that
$K'$ is not hyperbolic. One would expect that $S^3_{K'}(p'/q')$ is then non-hyperbolic for `generic' slopes $p'/q'$.
On the other hand, $S^3_{K}(p/q)$ is known to be hyperbolic for all but at most 10 slopes $p/q$, by work of the 
author and Meyerhoff [18]. Somewhat more pertinent to our analysis here is that $S^3_{K}(p/q)$ is hyperbolic 
when $|q| \geq 9$, a result also proved in [18]. 
However, it is in fact the case that
$S^3_{K'}(p'/q')$ might be hyperbolic for infinitely many slopes $p'/q'$, even without bound
on $|q'|$. A careful analysis of the JSJ decomposition for $S^3 - {\rm int}(N(K'))$ is
required before Theorem 1.2 can be proved in this case. This is completed in Section 5, after a
review of JSJ decompositions for knot exteriors is given in Section 4.

To prove Theorem 1.1, we must also establish that $S^3_K(p/q)$ `remembers' the Dehn surgery
when $K$ is not hyperbolic.
In the case where the JSJ piece of $S^3 - {\rm int}(N(K))$ containing $\partial N(K)$
is hyperbolic, the argument is again that this piece remains hyperbolic after Dehn surgery and the 
core of the attached solid torus is its shortest geodesic. But when the JSJ piece of $S^3 - {\rm int}(N(K))$
containing $\partial N(K)$ is non-hyperbolic, an alternative argument is required.
This piece is Seifert fibred, and we show that the Seifert fibration extends over the
attached solid torus. The core curve becomes an exceptional fibre and its singularity order
(which is the number of times that the regular fibres wind around it) is the largest in the
manifold. It is in this sense that the manifold `remembers' the Dehn surgery. In order
for this argument to work, we require here that $|p| \leq |q|$ and that $|q|$ is sufficiently large.
The proof is given in Section 6.

In this paper, all 3-manifolds will be oriented and connected. We use the symbol $\cong$ to denote the existence
of an orientation-preserving homeomorphism. Knots will not be oriented. We use
$(M,K)$ denote a pair consisting of the oriented 3-manifold $M$ and the unoriented
knot $K$ in $M$. Thus, $(M,K) \cong (M',K')$ means that there is an orientation-preserving
homeomorphism taking $M$ to $M'$ and taking $K$ to $K'$, but with no regard to an
orientation on $K$ or $K'$. We also use the terminology $(M,K) \cong_\partial (M,K')$
to denote the existence of an orientation-preserving homeomorphism of $M$
that takes $K$ to $K'$ and that equals the identity on $\partial M$.
More generally, if $M$ and $M'$ are 3-manifolds and a homeomorphism
between their boundaries has already been given, then we write $M \cong_\partial M'$
to denote the existence of a homeomorphism between $M$ and $M'$ that
restricts to the given homeomorphism between $\partial M$ and $\partial M'$.

If $\sigma$ is a slope on a toral boundary component of $M$, we let
$M(\sigma)$ denote the manifold that is obtained by Dehn filling along it.
If $K$ is a knot in $M$ and $\sigma$ is a slope on $\partial N(K)$, then the
result of performing Dehn surgery along $K$ using the slope $\sigma$ is
denoted $M_K(\sigma)$.

We say that a compact orientable 3-manifold is {\sl hyperbolic} if its interior admits a complete,
finite-volume hyperbolic structure. We say that a knot $K$ in a compact orientable 3-manifold $M$
is {\sl hyperbolic} if $M - {\rm int}(N(K))$ is hyperbolic.

\vskip 18pt
\centerline{\caps 2. Background results on Dehn surgery}
\vskip 6pt

The following is a version of Thurston's hyperbolic Dehn surgery theorem.

\noindent {\bf Theorem 2.1.} {\sl Let $X$ be a compact orientable hyperbolic 3-manifold
with boundary a non-empty collection of tori. Let $T_1, \dots, T_n$ be boundary components of $\partial X$. 
For each positive integer $i$, let $(\sigma_1^i, \dots, \sigma_n^i)$ be a collection of slopes, where
$\sigma_j^i$ lies on $T_j$. Suppose that $\sigma_j^i \not= \sigma_j^{i'}$ if $i \not= i'$. Then, for
each sufficiently large $i$, $X(\sigma^i_1, \dots, \sigma_n^i)$ is hyperbolic and the cores of the attached
solid tori are all geodesics. Furthermore, these geodesics all have lengths that tend to zero as $i \rightarrow \infty$,
whereas there is an $\epsilon > 0$ independent of $i$ such that all other geodesics in 
$X(\sigma^i_1, \dots, \sigma_n^i)$ that are not powers of a core curve have length at least $\epsilon$. Moreover, for any horoball
neighbourhood $N$ of the cusps of $X$, there is a subset $N_i$ of $X(\sigma^i_1, \dots, \sigma_n^i)$ consisting
of horoball neighbourhoods of the cusps and solid toral neighbourhoods of the surgery curves, 
so that the inclusion $X - N \rightarrow X(\sigma^i_1, \dots, \sigma_n^i) - N_i$ has 
bilipschitz constant that tends to $1$ as $i \rightarrow \infty$.}

The usual version of the hyperbolic Dehn surgery theorem ([30] or Theorem E.5.1 in [3]) asserts that there is a finite set of slopes
on each component of $\partial X$ so that, as long as these slopes are avoided when Dehn filling is performed,
the resulting manifold is hyperbolic. This is proved by perturbing the hyperbolic structure on $X$
to form an incomplete hyperbolic structure, and then taking the metric completion. If the perturbation is
chosen correctly, the points that are added in when forming the completion are a union of disjoint simple closed geodesics in
the resulting hyperbolic manifold. One can check that for most filling slopes, these geodesics are very short.
More precisely, if one chooses a basis for the first homology of each boundary torus so that
slopes on the torus can be identified with fractions $p/q \in {\Bbb Q} \cup \{ \infty \}$, then
the length of the corresponding geodesic is $O((|p|^2 + |q|^2)^{-1})$, according to an analysis by Neumann
and Zagier [22]. Thus, for all $i$ sufficiently large, $X(\sigma^i_1, \dots, \sigma_n^i)$ is hyperbolic and the cores of the attached
solid tori are all geodesics. Furthermore, these geodesics all have lengths that tend to zero as $i \rightarrow \infty$.
By the way that the hyperbolic structure on $X(\sigma_1^i, \dots, \sigma_n^i)$
is constructed, the hyperbolic manifolds $X(\sigma_1^i, \dots, \sigma_n^i)$ converge to
$X$ in the geometric topology (see Proposition E.6.29 in [3]). Hence, by Theorem E.2.4 in [3],
provided $\epsilon$ is small enough, the $\epsilon$-thick parts of $X$ and $X(\sigma^i_1, \dots, \sigma_n^i)$
are homeomorphic for all sufficiently large $i$. We deduce that the $\epsilon$-thin parts of
these manifolds have the same number of components. We choose $\epsilon$ to be smaller than the length of any
closed geodesic in $X$ and smaller than the 3-dimensional Margulis constant.
Therefore, the number of components of the $\epsilon$-thin part of $X$ is its
number of cusps. We can easily observe this many components of the $\epsilon$-thin
part of $X(\sigma^i_1, \dots, \sigma_n^i)$: its cusps and the solid toral neighbourhoods
of the surgery curves. Hence, this is precisely the $\epsilon$-thin part of $X(\sigma^i_1, \dots, \sigma_n^i)$.
Therefore, the only closed geodesics in $X(\sigma^i_1, \dots, \sigma_n^i)$
with length less than $\epsilon$ are powers of a core curve. Let $N$ be any horoball
neighbourhood of the cusps of $X$. This contains the $\delta$-thin part of $X$ for some
$\delta \leq \epsilon/2$. The proof of Theorem E.2.4 in [3]
provides a map $\phi_i$ from the $\delta$-thick part of $X$ to $X(\sigma^i_1, \dots, \sigma_n^i)$
that is a bilipschitz homeomorphism onto its image, with bilipschitz constant that tends
to $1$ as $i \rightarrow \infty$. The image of $\phi_i$ contains
the $2 \delta$-thick part of $X(\sigma^i_1, \dots, \sigma_n^i)$ and is homeomorphic to it. The complement
of $\phi_i(X-N)$ is not quite a horoball neighbourhood of the cusps together with solid toral neighbourhoods
of the core curves. But a minor adjustment to $\phi_i$ can be made to achieve this,
without modifying its bilipschitz constant too much. Hence, the final conclusion of the
theorem can be achieved. See [30, 3] for more details.

J{\o}rgensen and Thurston analysed the space of hyperbolic 3-manifolds with a bound on their volume.
Using the thick-thin decomposition of these manifolds, they were able to show the following (see [30] or Theorem E.4.8 in [3]).

\noindent {\bf Theorem 2.2.} {\sl Let $M_i$ be a sequence of 
pairwise distinct, oriented hyperbolic
$3$-manifolds with volume at most $V$. Then there is an oriented hyperbolic 3-manifold $X$ with volume at most $V$ and a collection
$T_1, \dots, T_n$ of toral boundary components of $\partial X$ with the following property.
For some subsequence of the $M_i$, there is a collection of slopes $(\sigma_1^i, \dots, \sigma_n^i)$, where
$\sigma_j^i$ lies on $T_j$, such that $M_i \cong X(\sigma_1^i, \dots, \sigma_n^i)$. Furthermore,
$\sigma_j^i \not= \sigma_j^{i'}$ if $i \not= i'$. }

In the theory of Dehn surgery, a key quantity that is considered is the {\sl distance} $\Delta(\sigma, \sigma')$
between two slopes $\sigma$ and $\sigma'$ on a torus. This is defined to be the modulus of the algebraic
intersection number between curves representing $\sigma$ and $\sigma'$.
A major theme is to produce upper bounds on $\Delta(\sigma, \sigma')$ when $\sigma$
and $\sigma'$ are slopes on the same boundary torus of a hyperbolic 3-manifold 
such that Dehn filling along these slopes results in non-hyperbolic 3-manifolds.

An important geometric quantity associated with any slope on the boundary of a hyperbolic
3-manifold is its `length' which is defined as follows. Let $M$ be a compact orientable 3-manifold with
interior that admits a complete, finite-volume
hyperbolic metric and let $N$ be a horoball neighbourhood of its cusps. Then the tori
$\partial N$ inherit a Euclidean metric. The {\sl length} of a slope $\sigma$ on $\partial M$
is the shortest curve on $\partial N$ with slope $\sigma$. Note that the length
depends on the choice of horoball neighbourhood $N$. It is usually advantageous
to take $N$ to be as large as possible, and so one often considers a {\sl maximal horoball neighbourhood}
of its cusps, which is by definition an open horoball neighbourhood that is not contained
in a larger one. Then $\partial N$ does not form a union of disjoint embedded tori. Instead,
these are immersed tori, but they still inherit a Euclidean metric, and so one can still
use $N$ to define the length of any slope on $\partial M$.

The relevance of slope length was first made apparent by work of Gromov and Thurston,
who proved that if $(\sigma_1, \dots, \sigma_n)$ is a collection of slopes on distinct boundary
components of $M$ and there is a horoball neighbourhood of the cusps with respect to which
each $\sigma_i$ has length more than $2 \pi$, then $M(\sigma_1, \dots, \sigma_n)$ admits
a negatively curved Riemannian metric [4]. As a consequence of Perelman's proof of the geometrisation
conjecture [25, 26, 27], such a manifold is now also known to be hyperbolic. The constant $2 \pi$ was
reduced to $6$ by the author [17] and Agol [1]. At the time, the conclusion on 
$M(\sigma_1, \dots, \sigma_n)$ was somewhat weaker, but again as a consequence
of geometrisation, $M(\sigma_1, \dots, \sigma_n)$ is now known to be hyperbolic if each
$\sigma_i$ has length more than $6$.

The following result relates slope length to the distance between slopes.

\vfill\eject
\noindent {\bf Theorem 2.3.} {\sl Let $M$ be a compact orientable non-hyperbolic 3-manifold, and let
$K$ be a hyperbolic knot in $M$. Let $N$ be a maximal horoball neighbourhood in $M - {\rm int}(N(K))$
of the cusp around $\partial N(K)$. Then, with respect to $N$, any slope $\sigma$ on $\partial N(K)$ satisfies
$${\rm length}(\sigma) \geq (\sqrt{3}/6) \ \Delta(\sigma, \mu),$$
where $\mu$ is the meridional slope on $\partial N(K)$. 
}

This was essentially proved in Corollary 2.4 of [7], but there the constant was $(\sqrt 3 /2 \pi)$ rather
than $(\sqrt 3 / 6)$. The improvement in the constant is due to the use of the 6-theorem rather
than the Gromov-Thurston $2\pi$ theorem.

The first inequality in the following theorem was proved by Futer, Kalfagianni and Purcell [8]. A similar result was proved by Cooper and the author [7]
but with less explicit estimates. The second inequality is due to Thurston [30].

\noindent {\bf Theorem 2.4.} {\sl Let $X$ be a compact orientable hyperbolic 3-manifold
with boundary a non-empty collection of tori. Let $T_1, \dots, T_n$ be toral boundary components
of $\partial X$ where $n > 0$. Let $(\sigma_1, \dots, \sigma_n)$ be a collection of slopes, where $\sigma_j$
lies on $T_j$. Suppose that there is a horoball neighbourhood of the cusps of $X$
on which each $\sigma_j$ has length at least $\ell_{\rm min} > 2\pi$. Then
$X(\sigma_1, \dots, \sigma_n)$ is a hyperbolic manifold satisfying
$$\left( 1 - \left ( {2\pi \over \ell_{\rm min}} \right )^2 \right)^{3/2} {\rm vol}(X) \leq
{\rm vol}(X(\sigma_1, \dots, \sigma_n)) < {\rm vol}(X).$$

}

The version of the above theorem in [7] was a crucial step in the proof of the following result of 
Cooper and the author (Theorem 4.1 in [7]).

\noindent {\bf Theorem 2.5.} {\sl Let $M$ be a compact orientable 3-manifold with boundary
a (possibly empty) union of tori. Let $\epsilon > 0$. Then there are only finitely many
compact orientable hyperbolic 3-manifolds $X$ and slopes $\sigma_1, \dots, \sigma_n$
on distinct components of $\partial X$ such that $M \cong X(\sigma_1, \dots, \sigma_n)$
and where the length of each $\sigma_i$ is at least $2 \pi + \epsilon$, when measured using
some horoball neighbourhood of the cusps of $X$.}

The following is due to the author and Meyerhoff [18]. 

\noindent {\bf Theorem 2.6.} {\sl Let $K$ be a hyperbolic knot in $S^3$. Then $S^3_K(p/q)$ is also
hyperbolic if $|q| \geq 9$.}

Theorem 1.2 in [18] states that 
if a 1-cusped hyperbolic 3-manifold is Dehn filled in two different ways and the resulting
manifolds are not hyperbolic, then the distance between the two surgery slopes is at most $8$.
Theorem 2.6 follows since $S^3_K(\infty) = S^3$ is not hyperbolic and $|q| = \Delta(p/q, \infty)$.

The following is due to Gordon and Luecke [12].

\noindent {\bf Theorem 2.7.} {\sl Let $K$ be a hyperbolic knot in $S^3$. Then $S^3_K(p/q)$ 
does not contain an incompressible embedded torus if $|q| > 2$.}

The following result generalises the above two theorems from the exterior of a knot in the 3-sphere
to knot exteriors in certain other 3-manifolds.

\noindent {\bf Theorem 2.8.} {\sl Let $M$ be the exterior of the unlink or unknot in $S^3$, and let $K$ be
a hyperbolic knot in $M$. Then $M_K(p/q)$ is irreducible and has incompressible boundary if $|q| \geq 2$.
Moreover, $M_K(p/q)$ is also hyperbolic provided $|q| \geq 3$.}

\noindent {\sl Proof.} The manifold $M_K(\infty) = M$ has compressible boundary. By Thurston's geometrisation
theorem [31], $M_K(p/q)$ is hyperbolic provided it contains no essential disc, sphere, annulus or torus, since its
boundary is a non-empty union of tori.
In the case where $M_K(p/q)$ has an essential disc, $|q| = \Delta(p/q, \infty) \leq 1$ by a theorem of Wu [33].
When $M_K(p/q)$ has an essential sphere, $|q| = \Delta(p/q, \infty) = 0$, by a theorem of Scharlemann [29].
If $M_K(p/q)$ has an essential annulus, then $|q| = \Delta(p/q, \infty) \leq 2$ by a result of Gordon and Wu [14].
When $M_K(p/q)$ has an essential torus, $|q| = \Delta(p/q, \infty) \leq 2$ by Gordon and Luecke [13]. $\square$

Gordon and Luecke [11] proved that knots in the 3-sphere are determined by their complements.
Specifically, any homeomorphism between $S^3 - {\rm int}(N(K))$ and $S^3 - {\rm int}(N(K'))$
for non-trivial knots $K$ and $K'$ in $S^3$ must send meridians to meridians. The following generalises
this to knot exteriors in some other 3-manifolds.

\noindent {\bf Theorem 2.9.} {\sl Let $M$ be the exterior of the unlink or unknot in $S^3$, and let $K$ and $K'$ be
knots in $M$ such that $M - {\rm int}(N(K))$ and $M - {\rm int}(N(K'))$ are irreducible. 
Suppose that there is a homeomorphism
between $M - {\rm int}(N(K))$ and $M - {\rm int}(N(K'))$ that is equal to the identity on $\partial M$.
Then this homeomorphism sends the meridian of $\partial N(K)$ to the meridian of $\partial N(K')$,
and therefore extends to a homeomorphism $(M,K) \cong_\partial (M, K')$.}

\noindent {\sl Proof.} The homeomorphism $M - {\rm int}(N(K')) \rightarrow M - {\rm int}(N(K))$ sends
the meridional slope on $\partial N(K')$ to some slope $\sigma$ on $\partial N(K)$. We want to show that
$\sigma$ is the meridional slope $\mu$ on $\partial N(K)$. 

When $M$ is the exterior of the unlink with at least two components, the deduction that $\sigma = \mu$
follows from a theorem of Scharlemann [29]. This is because $M_K(\mu) = M$ has compressible boundary,
whereas $M_K(\sigma) \cong M$ has a reducing sphere.

Suppose that $M$ is the exterior of the unknot. Suppose that $\sigma \not= \mu$. Then non-trivial surgery
on the knot $K$ in the solid torus $M$ yields the solid torus. Note also that $M - {\rm int}(N(K))$ is irreducible.
Gabai [9] showed that such knots $K$ have to be 1-bridge braids, and in particular must have non-zero winding number in $M$.
Now the homeomorphism $M_K(\sigma) \cong M$ is assumed to be equal to the identity on $\partial M$.
In particular, the slopes on $\partial M$ that bound discs in $M$ and $M_K(\sigma)$ are equal.
On the hand, for simple homological reasons, when non-trivial surgery along a knot in the solid torus $M$
with non-zero winding number is performed, the slopes on $\partial M$ that are
homologically trivial before and after the surgery are distinct. This is a contradiction,
and hence $\sigma = \mu$. $\square$

\vskip 18pt
\centerline {\caps 3. The case where both knots are hyperbolic}
\vskip 6pt

In this section, we prove the following generalisation of Theorem 1.3. This not only deals with knots in the 3-sphere,
but also knots in the exterior of the unknot or unlink. This extra level of generality will be useful when
handling non-hyperbolic knots in later sections.

\noindent {\bf Theorem 3.1.} {\sl Let $M$ be $S^3$ or the exterior of the unknot or unlink in $S^3$,
and let $K$ be a hyperbolic knot in $M$. Then there is a constant $C(K)$ with the following property. 
If $M_K(p/q) \cong_\partial M_{K'}(p'/q')$ for some hyperbolic knot $K'$ in $M$ and $|q'|\geq C(K)$, then 
$(M, K) \cong_\partial (M, K')$ and $p/q = p'/q'$.}

\noindent {\sl Proof.} Suppose that there is no constant $C(K)$ as in the theorem. Then there
is a sequence of slopes $p_i/q_i$ on $\partial N(K)$, a sequence of hyperbolic knots $K'_i$ in $M$ and slopes
$p'_i/q'_i$ on $\partial N(K'_i)$ such that 
\item{(i)} $|q_i'| \rightarrow \infty$;
\item{(ii)} $M_K(p_i/q_i) \cong_\partial M_{K'_i}(p'_i/q'_i)$; and
\item{(iii)} $(M,K) \not\cong_\partial (M,K'_i)$ or $p_i/q_i \not= p'_i/q'_i$.


Suppose first that the manifolds $M - {\rm int}(N(K'_i))$ run over infinitely many distinct homeomorphism types.
Then we may pass to a subsequence where they are pairwise distinct.
By Theorem 2.3, the length of $p'_i/q'_i$ tends to infinity and hence is larger
than $4 \pi$ for sufficiently large $i$. So by Theorem 2.4, the hyperbolic volume of
$M - {\rm int}(N(K'_i))$ is at most $(4/3)^{3/2}$ times that of $M_{K'_i}(p'_i/q'_i)$.
We are assuming that $M_{K'_i}(p'_i/q'_i)$ is homeomorphic to $M_K(p_i/q_i)$ and hence they have equal volume.
Since hyperbolic volume goes down when Dehn filling is performed, the volume of $M_K(p_i/q_i)$ is less than that
of $M - {\rm int}(N(K))$. So, if we set $V$ to be $(4/3)^{3/2}$ times the volume of $M - {\rm int}(N(K))$,
then the volume of $M - {\rm int}(N({K'_i}))$
is at most $V$. By Theorem 2.2, there is an oriented hyperbolic 3-manifold $X$ with volume at most $V$ and a collection
$T_1, \dots, T_n$ of toral boundary components of $\partial X$ with the following property.
For some subsequence of the $K'_i$, there is a collection of slopes $(\sigma_1^i, \dots, \sigma_n^i)$, where
$\sigma_j^i$ lies on $T_j$, such that $M - {\rm int}(N({K'_i})) \cong X(\sigma_1^i, \dots, \sigma_n^i)$. Furthermore,
$\sigma_j^i \not= \sigma_j^{i'}$ if $i \not= i'$. Hence, $M_{K'_i}(p'_i/q'_i) \cong 
X(\sigma_1^i, \dots, \sigma_n^i, \sigma^i)$ for some slope $\sigma^i$ on $\partial X - (T_1 \cup \dots \cup T_n)$. 

We now consider the possibility that the manifolds $M - {\rm int}(N(K'_i))$ run over only finitely many homeomorphism types.
We may then pass to a subsequence where they are a single oriented manifold $X$. So, in all cases,
$M_{K'_i}(p'_i/q'_i) \cong X(\sigma_1^i, \dots, \sigma_n^i, \sigma^i)$, 
for some slope $\sigma^i$ on $\partial X - (T_1 \cup \dots \cup T_n)$, but we allow the
possibility that $n = 0$.

By Theorem 2.1, there is a horoball neighbourhood $N$ of the cusps of $X$ and
subsets $N_i$ of $M - {\rm int}(N(K'_i))$ consisting of horoball neighbourhoods of the cusps
and solid toral neighbourhoods of the surgery curves,
so that the inclusion $X - N \rightarrow (M - {\rm int}(N(K'_i)) - N_i$ has bilipschitz constant that tends to $1$ as $i \rightarrow \infty$.
So the length of $\sigma^i$ as measured in $X$ tends to infinity as $i \rightarrow \infty$
because the length of $p'_i/q'_i$ tends to infinity.
We may therefore pass to a subsequence where the slopes $\sigma^i$ represent pairwise
distinct slopes on $\partial X$.
By Theorem 2.1, the core curves in $X(\sigma_1^i, \dots, \sigma_n^i, \sigma^i)$ are geodesics
when $i$ is sufficiently large,
with lengths that tend to zero. Similarly, applying Theorem 2.1 again, there is some $\epsilon > 0$ such that
the only geodesic in $M_K(p/q)$ with length less than $\epsilon$ is the core curve, provided
$|p| + |q|$ is greater than some constant $c$. There are only finitely many slopes with $|p| + |q| \leq c$,
and for each of these slopes $p/q$ that give hyperbolic fillings, we can compute the length of the shortest geodesic 
of the resulting manifold. If set $\epsilon$ to be less than this length, we deduce
that whenever $M_K(p_i/q_i)$ is hyperbolic, there is at most one geodesic in 
$M_K(p_i/q_i)$ with length less than $\epsilon$ and if there is such a geodesic, then it is the core curve.
Now, when $i$ sufficiently large, there are $n+1$ geodesics in $X(\sigma_1^i, \dots, \sigma_n^i, \sigma^i)$
with length less than $\epsilon$. This is a contradiction unless $n = 0$. Thus, $M - {\rm int}(N(K_i'))$ is homeomorphic to a fixed
oriented hyperbolic manifold $X$. 
The homeomorphism $M_{K_i'}(p_i'/q_i') \cong_\partial M_K(p_i/q_i)$ takes shortest geodesic to shortest geodesic
and so takes core curve to core curve.
So,  $M - {\rm int}(N({K'_i})) \cong M - {\rm int}(N(K))$ and the homeomorphism takes $p'_i/q'_i$ to $p_i/q_i$.
It restricts to the identity on $\partial M$.
It also takes the meridian of $M - {\rm int}(N(K'_i))$ to the meridian of $M - {\rm int}(N(K))$ by Theorem 2.9. 
So it extends to a orientation-preserving homeomorphism $(M, K'_i) \cong_\partial (M, K)$.
It takes the longitude on $\partial N(K'_i)$ to the longitude on $\partial N(K)$ for homological reasons. Hence, $p_i'/q'_i = p_i/q_i$. 
This contradicts our final hypothesis. $\square$

\vskip 12pt
\centerline {\caps 4. The JSJ decomposition of knot exteriors}
\vskip 6pt

In this section, we recall some details
about the JSJ decomposition for a knot exterior. A torus properly embedded
in a compact orientable irreducible 3-manifold $M$ is {\sl essential} if it is incompressible and not boundary parallel.
A properly embedded torus is a {\sl JSJ torus} if it is essential
and, moreover, any essential torus can be ambient isotoped off it.
The {\sl JSJ decomposition} for $M$ consists of one copy, up to
ambient isotopy, of each JSJ torus. The JSJ tori can be made simultaneously
disjoint, and so the JSJ decomposition is a union of disjoint essential tori.
When $M$ has boundary a (possibly empty) union of tori, the complement
of an open regular neighbourhood of the JSJ tori consists of 3-manifolds that are either Seifert fibred or
hyperbolic. Each component of this 3-manifold is called a {\sl JSJ piece}.

When the exterior of a knot $K$ in $S^3$ has no JSJ tori, it is hyperbolic or Seifert fibred.
The latter case arises exactly when $K$ is the unknot or a torus knot. Thus, we focus
on the situation where $S^3 - {\rm int}(N(K))$ has at least one JSJ torus. This was
analysed in detail by Budney [6]. The following summarises some of Theorem 4.18 in [6].

\noindent {\bf Theorem 4.1.} {\sl Let $K$ be a knot in $S^3$ such that $S^3 - {\rm int}(N(K))$ has at least one JSJ torus. 
Then any JSJ piece has one of the following forms:
\item{(1)} a Seifert fibre space with base space an annulus and with one exceptional fibre;
\item{(2)} a Seifert fibre space with base a planar surface with at least three boundary components and no exceptional fibres;
\item{(3)} a hyperbolic 3-manifold that is homeomorphic to the exterior of a hyperbolic link $L$ in $S^3$ such that when one specific component of $L$
is removed, the result is the unlink;
\item{(4)} the exterior of a torus knot in the 3-sphere.

}

The Seifert fibre spaces in Case (1) are called {\sl cable spaces}. One way in which they arise is when a knot $K$
is a cable of a knot $K'$, as follows. Let $K'$ be a knot in $S^3$, and let $N(K')$ be a solid torus regular neighbourhood
of $K'$. Then $K$ is a {\sl cable} of $K'$ if it lies on $\partial N(K')$ and the modulus of its winding number in $N(K')$ is
at least two. More precisely, we say that it is the {\sl $(r,s)$-cable} of $K$ if its slope on $\partial N(K')$ is $r/s$. The exterior
of $K$ then contains a cable space. For if we isotope $K$ a little into $N(K')$, then the result of removing a small open regular neighbourhood
of $K$ from $N(K')$ is a cable space. The core curve of $N(K')$ is the exceptional fibre, with singularity order $|s|$.

When there is a sequence of knots $K' = K_1, \dots, K_n = K$ such that each $K_{i+1}$ is a cable of $K_i$, we say that
$K$ is an {\sl iterated cable} of $K'$.

It is also worth recording the following extension of Theorem 4.1 when the JSJ piece contains $\partial N(K)$.

\noindent {\bf Theorem 4.2.} {\sl Let $K$ be a knot in $S^3$ such that $S^3 - {\rm int}(N(K))$ has at least one JSJ torus. 
Then the JSJ piece containing $\partial N(K)$ has one of the following forms:
\item{(1)} a Seifert fibre space with base space an annulus and with one exceptional fibre; in this case, $K$ is a cable of non-trivial knot;
\item{(2)} a Seifert fibre space with base a planar surface with at least three boundary components and no exceptional fibres; 
in this case, $K$ is a composite knot and each regular fibre on $\partial N(K)$ has meridional slope;
\item{(3)} a hyperbolic 3-manifold that is homeomorphic to the exterior of a hyperbolic link $L$ in $S^3$ such that when the component of $L$ corresponding
to $\partial N(K)$ is removed, the result is the unknot or unlink; the homeomorphism from the JSJ piece to $S^3 - {\rm int}(N(L))$
sends each slope $p/q$ on $\partial N(K)$ to the slope $p/q$ on the relevant component of $\partial N(L)$, and
it sends each slope $p/q$ on every other boundary torus to the slope $q/p$ on the relevant component of $\partial N(L)$.

}

\noindent {\sl Proof.} Most of this is contained in Theorem 4.18 of [6]. Case (1) of Theorem 4.1 corresponds to
conclusion (1)(b) of Theorem 4.18 in [6]. There, the statement is that $K$ is obtained by `splicing'
a non-trivial knot $K'$ with a `Seifert link' $S^{(p,q)}$. This is just the assertion that $K$ is the
$(q,p)$ cable of $K'$.

Case (2) of Theorem 4.1 corresponds to conclusion 1(c) of Theorem 4.18 in [6]. In this case,
$K$ is obtained by splicing some non-trivial knots with a `key-chain link'. (See Definition 4.8 in [6] 
for the definition of splicing.) This gives
that $K$ is a connected sum of these knots. The complement of the key-chain link
is the product of a planar surface $P$ and the circle $S^1$. The splicing construction gives
that the fibres of the Seifert fibration, which are of the form $\{ \ast \} \times S^1$, form
meridians on $\partial N(K)$.

The claims in Case (3) are essentially contained in [6]. When the component of $L$ corresponding
to $\partial N(K)$ is removed, the result is the unknot or unlink; this is contained in the definition of a `KGL' in [6].
The assertion about slopes in (3) follows from the definition of `splicing' in [6], as demonstrated in 
Example 4.9 in [6]. Note that we parametrise slopes on each boundary component of the JSJ
piece of $S^3 - {\rm int}(N(K))$ by observing that this torus separates $S^3$ into a solid
torus containing $K$ and the exterior of a non-trivial knot. On the boundary of the non-trivial
knot exterior, there is a canonical longitude and meridian, and hence a canonical fraction
$p/q \in {\Bbb Q} \cup \{ \infty \}$ for each slope.

Note that Case (4) of Theorem 4.1 does not
arise because we are assuming that $S^3 - {\rm int}(N(K))$ has at least one JSJ torus and hence is not
a torus knot exterior.  $\square$

\vskip 18pt
\centerline{\caps 5. Hyperbolic fillings on non-hyperbolic knots}
\vskip 6pt

It is a major theme in surgery theory that many properties of a 3-manifold with toral boundary are
preserved when a boundary component is Dehn filled along a `generic' slope. So one
might expect that if $K$ is a non-hyperbolic knot in $S^3$, then $S^3_K(p/q)$ is also non-hyperbolic
if $|q|$ is greater than some universal constant. Slightly surprisingly, this is not true, but we can analyse precisely
when this phenomenon arises.


\noindent {\bf Theorem 5.1.} {\sl Let $K$ be a non-hyperbolic knot in $S^3$. Suppose that
$S^3_{K}(p/q)$ is hyperbolic for $|q| \geq 2$. Then $K$ is an iterated cable of a 
hyperbolic knot $K'$ and $S^3_{K}(p/q) \cong S^3_{K'}(p'/q')$ for some slope $p'/q'$
where $|q| < |q'|$ and $|q|$ divides $|q'|$.}

The following result goes some way to proving Theorem 5.1.

\noindent {\bf Proposition 5.2.} {\sl Let $K$ be a non-hyperbolic knot in $S^3$ other than the unknot. Suppose that
$S^3_{K}(p/q)$ does not contain an incompressible torus and that $|q| \geq 2$. Then $K$ is a cable of a knot $K'$.}

\noindent {\sl Proof.} Let $X$ be the JSJ piece of $S^3 - {\rm int}(N(K))$
adjacent to $\partial N(K)$. If $X$ is all of $S^3 - {\rm int}(N(K))$, then $K$ is a torus knot, which is a cable of the unknot. So, we may assume that
$S^3 - {\rm int}(N(K))$ contains at least one JSJ torus.
The possibilities for $X$ are given in Theorem 4.2.

Suppose that $X$ is the exterior of a hyperbolic knot in $M$, where $M$ is the exterior
of the unknot or unlink in $S^3$. By Theorem 2.8, $X(p/q)$ is irreducible and has incompressible
boundary, since $|q| \geq 2$. Thus, $S^3_K(p/q)$ contains an incompressible torus, contrary to hypothesis.

Suppose that $X$ is homeomorphic to $P \times S^1$ for a planar surface $P$ with at least three boundary components, 
and that curves of the form $\{ \ast \} \times S^1$
on $\partial N(K)$ are meridians. So when we Dehn fill $\partial N(K)$ along a non-meridional slope,
the Seifert fibration on $P \times S^1$ extends over the attached solid torus.
Hence, the boundary tori of this filled-in piece are still incompressible, contrary to our hypothesis.

The only remaining possibility is that $X$ is Seifert fibred 
with base space an annulus and with one exceptional fibre. Then,
$K$ is a cable of a knot $K'$. $\square$

Then an analysis of surgery along cabled knots gives the following.

\noindent {\bf Proposition 5.3.} {\sl Let $K$ be a cable of a non-trivial knot $K'$ in $S^3$. Suppose that 
$S^3_{K}(p/q)$ does not have an incompressible torus and that $|q| \geq 2$. Then 
$S^3_{K}(p/q) \cong S^3_{K'}(p'/q')$ for some slope $p'/q'$, 
where $|q| < |q'|$ and $|q|$ divides $|q'|$.}

\noindent {\sl Proof.} Since $K$ is a cable of $K'$, by definition $K$ lies on the torus $\partial N(K')$. 
Say that it is a curve with slope $r/s$ on $\partial N(K')$, where $|s| > 1$. As in the previous proof, this means that there is a cable 
space $X$ embedded in the exterior of $K$, with one boundary component being $\partial N(K)$. 
This cable space is a Seifert fibre space with base space an annulus, and with one exceptional fibre.
The exceptional fibre is the core of the solid torus bounded by $T = \partial X - \partial N(K)$.
The regular fibres are parallel to the slope $r/s$ on $\partial N(K')$. They have slope $rs$ on
$\partial N(K)$. When we Dehn fill $\partial N(K)$ along the slope $p/q$, we attach a solid
torus to $X$, forming a 3-manifold $Y$. The Seifert fibration on 
$X$ extends over the solid torus, as long as the filling slope is not equal to the slope of the
regular fibres. But we are assuming that $|q| > 1$, whereas the slope of the regular fibres
is the integer $rs$. Thus, $Y$ is Seifert fibred. Its base space is a disc, and it has either
one or two exceptional fibres. When it has two, $T$ is incompressible in $Y$. Hence,
$S^3_{K}(p/q)$ contains an essential torus, which is contrary to one of our assumptions.
Therefore, we deduce that when the solid torus is attached to $X$, the Seifert fibration on $X$
extends over this solid torus without introducing an exceptional fibre. Therefore, $X(p/q)$ is a
solid torus, and therefore $S^3_K(p/q)$ is homeomorphic to the manifold that is obtained
by Dehn filling $K'$ along some slope $p'/q'$.

We now determine the required properties of $q'$. Since the solid torus attached to $X$ is
fibred using only regular fibres, the distance between the surgery slope
and the slope of the regular fibres on $\partial N(K)$ is $1$. Now the slope of the regular fibres is
$rs$ and so the surgery slope $p/q$ must be of the form $(1+krs)/k$, for some integer $k$.
To see this, note that such slopes are precisely those that are, homologically, equal to
a meridian plus some multiple of the slope of the regular fibres.
Now, the Seifert fibration on $Y$ can be viewed as starting from $P \times S^1$,
where $P$ is a three-times punctured sphere, with fibres of the form $\{ \ast \} \times S^1$,
with two solid tori attached. We can perform an automorphism to $P \times S^1$
by performing a Dehn twist along an annulus, with one boundary component on
$\partial N(K)$ and the other boundary component on $T$. 
Applying such Dehn twists, we can take the slope $(1+krs)/k$ on $\partial N(K)$ back to
the meridian $1/0$, but it changes the way that $X$ is attached to $S^3 - {\rm int}(N(K'))$. The way that this changes is by applying $k$ Dehn twists
along the slope of the regular fibre. On $T$, this slope is $r/s$. Hence, the meridian
is taken to $(1+kr)/ks$. We therefore deduce that $S^3_{K}(p/q) = S^3_{K}((1+krs)/k)
\cong S^3_{K'}((1+kr)/ks)$. Setting $k = q$ and $ks = q'$, we obtain the proposition.
$\square$

\noindent {\sl Proof of Theorem 5.1.} We will prove this by induction on the number of JSJ
tori in $S^3 - {\rm int}(N(K))$. The induction starts where this number is zero.
Since $K$ is not hyperbolic, $S^3 - {\rm int}(N(K))$ must therefore be Seifert fibred. Hence, $K$ is a 
torus knot or the unknot. Every surgery on a torus knot or the unknot yields a non-hyperbolic 3-manifold,
and this contradicts our assumption that $S^3_{K}(p/q)$ is hyperbolic.

We now prove the inductive step. Let $K$ be a non-hyperbolic knot in $S^3$ other than the unknot. Suppose that
$S^3_{K}(p/q)$ is hyperbolic for $|q| \geq 2$. By Proposition 5.2, $K$ is a cable
of a knot $K'$. Note that $K'$ cannot be the unknot, because then $K$ would be a
torus knot and hence $S^3_K(p/q)$ could not be hyperbolic.
By Proposition 5.3 and the assumption that $S^3_{K}(p/q)$ is hyperbolic,
$S^3_{K}(p/q) \cong S^3_{K'}(p'/q')$ for some slope $p'/q'$,
where $|q| < |q'|$ and $|q|$ divides $|q'|$.
If $K'$ is hyperbolic, the theorem is proved. Therefore, suppose that $K'$
is not hyperbolic. By induction, $K'$ is an iterated cable of a hyperbolic knot $K''$
and $S^3_{K'}(p'/q') \cong S^3_{K''}(p''/q'')$ where $|q'| < |q''|$ and $|q'|$ divides $|q''|$.
So, $K$ is an iterated cable of $K''$ and $S^3_{K}(p/q) \cong S^3_{K''}(p''/q'')$
where $|q| < |q''|$ and $|q|$ divides $|q''|$.
$\square$

We are now in a position to prove Theorem 1.2.

\noindent {\sl Proof of Theorem 1.2.} Let $K$ be a hyperbolic knot. Let $C(K)$ be the constant
from Theorem 1.3. We show that any slope $p/q$ is characterising for $K$, provided that
$|q| \geq \max \{ C(K), 9 \}$.

Let $K'$ be a knot in $S^3$ such that $S^3_K(p/q) \cong S^3_{K'}(p/q)$. If $K'$ is hyperbolic,
then by Theorem 1.3, $(S^3,K) \cong (S^3, K')$ and the theorem is proved. So, suppose $K'$
is not hyperbolic. By Theorem 2.6, $S^3_K(p/q)$ is hyperbolic 
since $|q| \geq 9$. By Theorem 5.1, $K'$ is an iterated cable of a hyperbolic knot $K''$ and $S^3_{K'}(p/q) \cong 
S^3_{K''}(p''/q'')$ where $|q| < |q''|$ and $|q|$ divides $|q''|$. Since $|q''| \geq C(K)$,
we deduce from $S^3_K(p/q) \cong S^3_{K''}(p''/q'')$
that $(S^3, K) \cong (S^3, K'')$ and $p/q  = p''/q''$, using Theorem 1.3. In particular, $|q| = |q''|$, which
contradicts our earlier conclusion that $|q| < |q''|$. $\square$

\vskip 18pt
\centerline {\caps 6. Proof of the main theorem}
\vskip 6pt

In this section, we complete the proof of Theorem 1.1. Since Theorem 1.2 was proved in the last section, we
must focus on the case where $K$ is not hyperbolic. Our first result deals with the relationship between
the JSJ decompositions of $S^3 - {\rm int}(N(K))$ and $S^3_K(p/q)$ for suitable slopes $p/q$.

\noindent {\bf Theorem 6.1.} {\sl Let $K$ be a knot in $S^3$. Suppose that $|q| \geq 9$. In the case where
$K$ is a torus knot or a cable knot, suppose also that $|p| \leq |q|$. Then $S^3_K(p/q)$ is irreducible
and the JSJ tori for $S^3 - {\rm int}(N(K))$
form the JSJ tori for $S^3_K(p/q)$. Moreover, if the JSJ piece of $S^3 - {\rm int}(N(K))$ containing $\partial N(K)$ is hyperbolic, then so is the
corresponding piece in $S^3_K(p/q)$. If the JSJ piece of $S^3 - {\rm int}(N(K))$ containing $\partial N(K)$ is Seifert fibred, then so is the
corresponding piece in $S^3_K(p/q)$ and the core of the attached solid torus is an exceptional fibre with singularity order at least $|q|$.

}

\noindent {\sl Proof.} Let $Y$ be the JSJ piece for $S^3 - {\rm int}(N(K))$ that contains $\partial N(K)$.
We will show that if $Y$ is hyperbolic, then $Y(p/q)$ is too. We will also show that if $Y$ is Seifert fibred,
then this extends to a Seifert fibration for $Y(p/q)$ in which the core of the surgery solid torus
is an exceptional fibre with singularity order at least $|q|$. 
In both cases, we deduce that $Y(p/q)$ has incompressible boundary
and is not a copy of $T^2 \times I$. This will prove the theorem. This is
because the JSJ tori for $S^3 - {\rm int}(N(K))$ form a decomposition of $S^3_K(p/q)$ into pieces
with incompressible boundary, that are either hyperbolic or Seifert fibred and not $T^2 \times I$,
and when two Seifert fibred pieces meet along a torus, then their Seifert fibrations are not isotopic on this torus.
This is enough to deduce that these tori are the JSJ tori for $S^3_K(p/q)$ and that $S^3_K(p/q)$ is
irreducible.

Suppose first that $Y$ is hyperbolic. Then by Theorem 4.2, $Y$ is homeomorphic to
the exterior of a hyperbolic link $L$ in $S^3$, with the property that when the component of $L$
corresponding to $\partial N(K)$ is removed, the result is the unlink, the unknot or the empty set. 
Furthermore, the homeomorphism
$Y \rightarrow S^3 - {\rm int}(N(L))$ takes the slope $p/q$ on $\partial N(K)$
to the slope $p/q$ on the relevant component of $\partial N(L)$. By Theorem 2.6 or 2.8,
$Y(p/q)$ is hyperbolic if $|q| \geq 9$.

Suppose now that $Y$ is Seifert fibred. If $Y$ has more than two boundary components,
then by Theorem 4.2, it is of the form $P \times S^1$ for a planar surface $P$ and the slope of the fibres on $\partial N(K)$
is meridional. Hence, if $p/q$ surgery is performed and $q \not = 0$, then the
Seifert fibration extends over the solid torus. Moreover, if $|q| \geq 2$, then
the core of the solid torus becomes an exceptional fibre with order $|q|$, as required.

Consider now the case when $Y$ is Seifert fibred with two boundary components.
Then $Y$ is a cable space and $K$ is an $(r,s)$-cable of a non-trivial knot.
The slope of the regular fibre is $rs$ on $\partial N(K)$. The distance between it
and $p/q$ is $|qrs - p|$. Note that
$$|qrs - p| = |q| |rs - (p/q)| \geq |q|,$$
where the inequality follows from the fact that $|rs| \geq 2$ and $|p/q| \leq 1$.
We are assuming that $|q| \geq 2$. Hence, we deduce that the Seifert fibration
of $Y$ extends to $Y(p/q)$, and 
the core of the attached solid torus becomes an exceptional fibre with singularity
order at least $|q|$, again as required.

Suppose now that $Y$ is Seifert fibred with a single boundary component.
Then $Y$ is the exterior of an $(r,s)$ torus knot. Its Seifert fibration extends over the attached solid torus,
as long as the surgery slope is not equal to that of the regular fibre, which is
the integer $rs$. As in the case of the cable space, the distance between $p/q$ and $rs$
is at least $|q|$ and hence again the core of the attached solid torus becomes an exceptional fibre with singularity
order at least $|q|$.  $\square$

\noindent {\sl Proof of Theorem 1.1.} Let $K$ be a knot in $S^3$. If $K$ is the unknot,
then Theorem 1.1 follows from the main theorem of [16] by 
Kronheimer, Mrowka, Ozsv\'ath and Szab\'o. If $K$ is hyperbolic,
then Theorem 1.1 follows from Theorem 1.2, which we proved in the previous section.
If $K$ is a torus knot, then it is a theorem of
McCoy [20] that every non-integral slope, with at most finitely many exceptions, is characterising for $K$. 
So suppose that $S^3 - {\rm int}(N(K))$ contains at least one JSJ torus,
and let $Y$ be the JSJ piece containing $\partial N(K)$. Let $K'$ be another knot in $S^3$
and let $Y'$ be the JSJ piece of $S^3 - {\rm int}(N(K'))$ containing $\partial N(K')$.

Suppose that $S^3_{K'}(p/q) \cong S^3_K(p/q)$ for some slope $p/q$ where $|p| \leq |q|$ and 
$|q|$ is larger than some constant that depends on $K$, and which is to be determined.
We will initially assume that $|q| \geq 9$ so that Theorem 6.1 applies. So, $S^3_K(p/q)$
is irreducible and the JSJ tori for $S^3 - {\rm int}(N(K))$ form the JSJ tori for $S^3_K(p/q)$.
Theorem 6.1 also applies to $K'$
and hence the JSJ tori for $S^3 - {\rm int}(N(K'))$ form the JSJ tori for $S^3_{K'}(p/q)$.
The homeomorphism $S^3_{K'}(p/q) \cong S^3_{K}(p/q)$ takes JSJ tori to JSJ tori, after possibly
applying an isotopy. (This isotopy may move the surgery curves in the manifold.)

\noindent {\sl Case 1.} $Y'$ is Seifert fibred.

If $|q|$ is sufficiently large, then by Theorem 6.1, $Y'(p/q)$ is Seifert fibred and the core of the surgery solid torus
in $Y'(p/q)$ is an exceptional fibre with singularity order that 
is greater than that of any exceptional fibre in any Seifert fibred JSJ piece of
$S^3 - {\rm int}(N(K))$.
Hence, the homeomorphism $S^3_{K'}(p/q) \cong S^3_K(p/q)$ must take
this exceptional fibre to the core of the surgery solid torus in $S^3_K(p/q)$.
Thus, it restricts to $S^3 - {\rm int}(N(K')) \cong S^3 - {\rm int}(N(K))$ taking 
the slope $p/q$ to $p/q$. It also takes meridians to meridians, by the theorem
of Gordon and Luecke [11]. Hence, it extends to $(S^3, K') \cong (S^3, K)$, as required.

\noindent {\sl Case 2.} $Y'$ is hyperbolic.

By Theorem 4.2, there is a homeomorphism $Y' \rightarrow S^3 - {\rm int}(N(L'))$
for some hyperbolic link $L'$. Furthermore, if the component of $L'$ corresponding to
$\partial N(K')$ is removed, the result is the unknot or unlink. The homeomorphism
sends the slope $p/q$ on $\partial N(K')$ to the slope $p/q$ on the relevant component
of $\partial N(L')$.

We claim that the homeomorphism $S^3_{K'}(p/q) \cong S^3_{K}(p/q)$ takes the JSJ piece in $S^3_{K'}(p/q)$ containing the surgery curve
to the JSJ piece in $S^3_{K}(p/q)$ containing the surgery curve, provided that $|q|$
is greater than some constant that depends on $K$.
If not, then $Y'(p/q)$ is homeomorphic to a JSJ piece $X$ for $S^3 - {\rm int}(N(K))$
not containing $\partial N(K)$. Now, $S^3 - {\rm int}(N(K))$ has
only finitely many JSJ pieces. By Theorem 2.5, for each such piece $X$ and each $\epsilon > 0$, there are only finitely
many hyperbolic 3-manifolds $Y'$ and slopes $\sigma$ on $Y'$ with length at least $2 \pi + \epsilon$
such that $Y'(\sigma) \cong X$. Now the length of $p/q$ on $Y'$ is at least $(\sqrt{3}/6) |q| > 2 \pi$ if $|q| \geq 22$, by Theorem 2.3. 
Thus, by choosing $p/q$ to avoid finitely many values and so that $|q| \geq 22$, we can arrange
that $Y'(p/q) \not \cong X$ for any JSJ piece $X$ of $S^3 - {\rm int}(N(K))$. 
This proves the claim.

Thus the homeomorphism 
$S^3_K(p/q) \cong S^3_{K'}(p/q)$ restricts to a homeomorphism $Y(p/q) \cong Y'(p/q)$.
It also restricts to a homeomorphism $h \colon S^3 - N(K) - {\rm int} (Y) \cong S^3 - N(K') - {\rm int}(Y')$.
Let $Z$ and $Z'$ denote these manifolds.

Note that each boundary component of $Z$ is a torus in $S^3$,
which therefore bounds a solid torus on at least one side. The solid torus must contain
$K$, since its boundary torus is incompressible in $S^3 - {\rm int}(N(K))$. The same argument
gives that, on the other side of the torus, there is the exterior of a non-trivial knot.
Thus, each component of $Z$ is homeomorphic to
the exterior of a non-trivial knot in $S^3$. The same is true of each component of
$Z'$. Hence, the homeomorphism 
$h \colon Z \rightarrow Z'$
must send meridians to meridians by [11]. It must also send longitudes to
longitudes for homological reasons. Since it is orientation-preserving, it must therefore
send each slope $r/s$ on each boundary component to the slope $r/s$
on the image boundary component. We say that such a homeomorphism is
{\sl slope-preserving}. This restricts to a homeomorphism $h|\partial Z$ from the boundary
of $Y(p/q)$ to the boundary of $Y'(p/q)$.

Now, $Y'(p/q)$ is hyperbolic by Theorem 6.1. Thus, $Y(p/q)$ is also, and hence so is $Y$ by Theorem 6.1.
So by Theorem 4.2, $Y$ is homeomorphic to the exterior of a
hyperbolic knot $K_\ast$ in $M$, where $M$ is the exterior of the unknot or unlink
in $S^3$. The homeomorphism takes $\partial N(K)$ to $\partial N(K_\ast)$ and sends each slope $r/s$ on $\partial N(K)$
to $r/s$ on $\partial N(K_\ast)$. On the remaining components of $\partial Y$, the homeomorphism sends each
slope $r/s$ to the slope $s/r$ on the relevant component of $\partial M$.
Similarly, $Y'$ is homeomorphic to the exterior of a hyperbolic knot $K_\ast'$
in $M$, and the homeomorphism has the same effect on slopes as was the case
with $Y$. We obtain an induced homeomorphism $\partial M \rightarrow \partial M$
as follows. We have a homeomorphism $M - {\rm int}(N(K_\ast)) \rightarrow Y$
which restricts to a homeomorphism $\partial M \rightarrow \partial Y - \partial N(K)$.
There is a homeomorphism $Y(p/q) \rightarrow Y'(p/q)$ which restricts to a
homeomorphism $\partial Y - \partial N(K) \rightarrow \partial Y' - \partial N(K')$.
And then there is a homeomorphism $\partial Y' - \partial N(K') \rightarrow \partial M$.
We observe that this is a slope-preserving homeomorphism $\partial M \rightarrow \partial M$.
We can ensure that this is actually equal to the identity on $\partial M$ by possibly changing the knot $K'_\ast$, as follows.
There is a homeomorphism $M \rightarrow M$ that achieves any given permutation 
of its boundary components, since one may permute the components of an unlink in
any way by an isotopy of the 3-sphere. There is also a homeomorphism $M \rightarrow M$
which acts as $-{\rm id}$ on one component of $\partial M$ and acts as the identity
on the remaining components. This is achieved by an isotopy of the 3-sphere, which moves
one component of the unlink back to itself but reversing its orientation.
By applying such homeomorphisms of $M$ to itself, we take $K_\ast'$ to another
knot $K_\ast''$, but now the relevant homeomorphism $\partial M \rightarrow \partial M$
is the identity. In other words, $M_{K_\ast}(p/q) \cong_\partial M_{K_\ast''}(p/q)$.
So, by Theorem 3.1,
if $|q|$ is sufficiently large, $(M,K_\ast) \cong_\partial (M,K_\ast'')$. 
This gives a homeomorphism $S^3 - {\rm int}(Z) \rightarrow S^3 - {\rm int}(Z')$ taking
$K$ to $K'$, and that equals $h$ on its boundary.
This extends to a homeomorphism $(S^3,K) \cong (S^3,K')$ using the homeomorphism $h \colon Z \rightarrow Z'$.
$\square$

We note that the same argument gives the somewhat stronger, albeit less catchy, version of Theorem 1.1.

\noindent {\bf Theorem 6.2.} {\sl Let $K$ be a knot in $S^3$. Then the slope $p/q$ is characterising for $K$
provided that 
$$|q| \ \left ( \min \bigl \{ |(p/q) - n| : n \in {\Bbb Z} \backslash \{ -1, 0, 1 \} \bigr \} \right ) $$
is sufficiently large.

}

In Theorem 1.1, we required that $|p/q| \leq 1$, whereas in Theorem 6.2, we only require that,
in some sense, $p/q$ is not too close to any integer in ${\Bbb Z} \backslash \{ -1, 0, 1 \} $.

There is only one step in the proof where the hypothesis that $|p| \leq |q|$ is
used. It is in the proof of Theorem 6.1, when the JSJ piece of $S^3 - {\rm int}(N(K))$ containing $\partial N(K)$ is Seifert fibred.
The conclusion there is that this Seifert fibration extends over the surgery solid torus, with the core
being a singular fibre with singularity order at least $|q|$.
However, all one really needs is that this singularity order is at least some function $f(p/q)$
where $f$ is independent of $K$ and that, for certain values of $p/q$, this function $f(p/q)$ is 
large. Examining the proof of Theorem 6.1, we see that the singularity order is
$|qrs - p|$ for suitable coprime integers $r$ and $s$ where $|rs| \geq 2$.
Note that
$$|qrs - p| =  |q| |rs - (p/q)| \geq 
|q| \ \min \bigl \{ |(p/q) - n| : n \in {\Bbb Z} \backslash \{ -1, 0, 1 \} \bigr \}.$$
Hence if the latter quantity is large, then so is the order of the singular fibre.
The proof of Theorem 1.1 can then be readily adapted to establish Theorem 6.2. $\square$

\vskip 18pt
\centerline {\caps References}
\vskip 6pt


\item{1.} {\caps I. Agol,} {\sl Bounds on exceptional Dehn filling.} Geom. Topol. 4 (2000), 431--449. 


\item{2.} {\caps K. Baker, K. Motegi,} {\sl Non-characterizing slopes for hyperbolic knots},
arXiv:1601.01985 (2016).

\item{3.} {\caps R. Benedetti, C. Petronio,} {\sl Lectures on hyperbolic geometry.}
Universitext. Springer-Verlag, Berlin, 1992.

\item{4.} {\caps S. Bleiler, C. Hodgson,} {\sl Spherical space forms and Dehn filling.} Topology 35 (1996), no. 3, 809--833.

\item{5.} {\caps S. Bleiler, C. Hodgson, J. Weeks,} {\sl Cosmetic surgery on knots.}
Proceedings of the Kirbyfest (Berkeley, CA, 1998), 23Ð34, Geom. Topol. Monogr., 2, Geom. Topol. Publ., Coventry, 1999

\item{6.} {\caps R. Budney}, {\sl JSJ-decompositions of knot and link complements in $S^3$.} 
Enseign. Math. (2) 52 (2006) 319--359.


\item{7.} {\caps D. Cooper, M. Lackenby,} {\sl Dehn surgery and negatively curved 3-manifolds,} 
J. Differential Geom. 50 (1998) 591--624.

\item{8.} {\caps D. Futer, E. Kalfagianni, J. Purcell,} {\sl Dehn filling, volume, and the Jones polynomial.}
J. Differential Geom. 78 (2008), no. 3, 429--464.

\item{9.} {\caps D. Gabai,} {\sl Surgery on knots in solid tori.} Topology 28 (1989), no. 1, 1--6. 

\item{10.} {\caps C. Gordon,} {\sl Some aspects of classical knot theory,} Lecture Notes in Math. 685,
1--60, Springer-Verlag, New York, 1978.

\item{11.} {\caps C. Gordon, J. Luecke,} {\sl Knots are determined by their complements},
J. Amer. Math. Soc. 2 (1989), no. 2, 371--415.

\item{12.} {\caps C. Gordon, J. Luecke,} {\sl Dehn surgeries on knots creating essential tori, I,}
Comm. Anal. Geom. 3 (1995) 597-644.

\item{13.} {\caps C. Gordon, J. Luecke,} {\sl Toroidal and boundary-reducing Dehn fillings,} 
Topology Appl. 93 (1999) 77--90.

\item{14.} {\caps C. Gordon, Y.-Q. Wu,} {\sl Annular and boundary reducing Dehn fillings,} Topology 39 (2000) 531--548.

\item{15.} {\caps R. Kirby,} {\sl A calculus for framed links in $S^3$.} Invent. Math. 45 (1978), no. 1, 35--56.

\item{16.} {\caps P. Kronheimer, T. Mrowka, P. Ozsv\'ath, Z. Szab\'o,} {\sl 
Monopoles and lens space surgeries,} Ann.of Math. (2), 165 (2007) 457--546.

\item{17.} {\caps M. Lackenby}, {\sl Word hyperbolic Dehn surgery.} Invent. Math. 140 (2000), no. 2, 243--282.

\item{18.} {\caps M. Lackenby, R. Meyerhoff,} {\sl The maximal number of exceptional Dehn surgeries,} 
Invent. Math. 191 (2013) 341--382.

\item{19.} {\caps W. B. R. Lickorish,} {\sl A representation of orientable combinatorial 3-manifolds.}
Ann. of Math. (2) 76 (1962) 531--540.

\item{20.} {\caps D. McCoy}, {\sl Non-integer characterising slopes for torus knots}, arXiv:1610.03283 (2016).

\item{21.} {\caps G. Mostow,} {\sl Strong rigidity of locally symmetric spaces.}
Annals of Mathematics Studies, No. 78. Princeton University Press, Princeton, N.J.; University of Tokyo Press, Tokyo, 1973.

\item{22.} {\caps W. Neumann, D. Zagier,} {\sl Volumes of hyperbolic three-manifolds.} Topology 24 (1985), no. 3, 307--332. 

\item{23.} {\caps Y. Ni, X. Zhang,} {\sl Characterizing slopes for torus knots,} Algebr. Geom. Topol.
14 (2014) 1249--1274.

\item{24.} {\caps P. Ozsv\'ath, Z. Szab\'o,} {\sl The Dehn surgery characterization of the trefoil and the
figure eight knot,} arXiv:0604079 (2006).

\item{25.} {\caps G. Perelman,} {\sl The entropy formula for the 
Ricci flow and its geometric applications,} Preprint,
arxiv:math.DG/0211159

\item{26.} {\caps G. Perelman,} {\sl  Ricci flow with surgery on three-manifolds,}
Preprint, arxiv:math.DG/0303109

\item{27.} {\caps G. Perelman,} {\sl Finite extinction time for the 
solutions to the Ricci flow on certain three-manifolds,} Preprint,
arxiv:math.DG/0307245 

\item{28.} {\caps  Y. Rieck, Y. Yamashita,} {\sl Cosmetic surgery and the link volume of hyperbolic 3-manifolds.}
Algebr. Geom. Topol. 16 (2016), no. 6, 3445--3521

\item{29.} {\caps M. Scharlemann,} {\sl Producing reducible 3-manifolds by surgery on a knot,}
Topology 29 (1990) 481--500.

\item{30.} {\caps W. Thurston}, {\sl The geometry and topology of three-manifolds}, Notes (1980), available for download from
http://library.msri.org/books/gt3m/

\item{31.} {\caps W. Thurston}, {\sl Three-dimensional manifolds, Kleinian groups and hyperbolic geometry}. 
Bull. Amer. Math. Soc. (N.S.). 6 (1982) 357--381.

\item{32.} {\caps A. Wallace,} {\sl Modifications and cobounding manifolds,} Canad. J. Math. 12 (1960) 503--528.

\item{33.} {\caps Y.-Q. Wu,} {\sl Incompressibility of surfaces in surgered 3-manifolds,} Topology 31 (1992) 271--279.

\vskip 12pt
\+ Mathematical Institute, University of Oxford, \cr
\+ Radcliffe Observatory Quarter, Woodstock Road, Oxford OX2 6GG, United Kingdom. \cr

\end